\def\versiondate{24 May 1999}
\input math.macros
\input Ref.macros

\checkdefinedreferencetrue
\continuousfigurenumberingtrue
\theoremcountingtrue
\sectionnumberstrue
\forwardreferencetrue
\initialeqmacro

\def\st{:\,}
\def\ev#1{{\cal #1}} 
\def\A{\ev A}
\def\verts{{V}}
\def\edges{{E}}



\def\genassumps{Let $(\Xi, \F, \P, \Gamma)$ be a measure-preserving
   dynamical system
   and $p$ be a $\Gamma$-invariant random environment from $\Xi$}
\def\SS{{\cal S}}   

\def\Ph{\widehat \P}   

\def\F{{\cal F}}



\def\traj{\theta}

\def\inv{{\cal I}}    
\def\bdend{\zeta}     
\def\scnry{\Upsilon} 

\def\BLPSgip{Benjamini, Lyons, Peres, and Schramm (1999)}

\def\BLSpert{Benjamini, Lyons, and Schramm (1998)}
\def\LSindist{Lyons and Schramm (1998)}

\ifproofmode \relax \else\head{} {Version of \versiondate}\fi 
\vglue20pt

\title{Stationary Measures for Random Walks}
\title{in a Random Environment with Random Scenery}

\author{Russell Lyons and Oded Schramm}

\abstract{Let $\Gamma$ act on a countable set $\verts$ with only finitely
many orbits. Given a $\Gamma$-invariant random environment for a Markov
chain on $\verts$ and a random scenery, we exhibit, under certain
conditions, an equivalent stationary measure for the environment and scenery
from the viewpoint of the random walker. Such theorems have been very
useful in investigations of percolation on quasi-transitive graphs.
}

\bottom{Primary 
60B99, 
60J15.  
Secondary
28D15. 
} 
{Cayley graph, group, transitive. 
}
{Research partially supported by NSF grant DMS-9802663 (Lyons) and
the Sam and Ayala Zacks Professorial Chair (Schramm).}

\bsection {Introduction}{s.intro}

Given a state space for a Markov chain, one might assign transition
probabilities randomly in order to finish specifying the Markov chain. In
such a case, one speaks about {\bf random walk in a random environment}, or
{\bf RWRE} for short. If we do not condition on the transition probabilities, 
such a stochastic process is usually no longer a
Markov chain. The first investigation of RWREs is due to Solomon (1975).
Their properties are often surprising. 

Alternatively, given a completely specified Markov chain, which we shall
refer to as a random walk, there might be
a random field on the state space, i.e., a collection of random variables
indexed by the state space. This random field is called a {\bf random
scenery}. As the
random walker moves, he observes the scenery at his location. Perhaps the
first explicit investigation of random walks in random scenery was 
Lang and Nguyen (1983).

Of course, one may combine these processes to obtain a random walk in a
random environment with random scenery, or {\bf RWRERS} for short. This has
not been looked at much except in the case where the scenery arises from
percolation on a graph and determines the environment (H\"aggstr\"om
(1997), H\"aggstr\"om and Peres (1999), Lyons and Schramm (1998)). In fact,
the purpose of those investigations was to find out information about the
scenery; the corresponding RWRE was used as a tool to probe the scenery. 

In general, one would like a stationary probability measure on the
trajectories of an RWRERS that is equivalent to
(mutually absolutely continuous with)
the natural probability measure giving the environment, the
scenery, and the trajectory of the Markov chain given the environment.
Here, stationarity
means that when looked at from the viewpoint of the random walker, one
should see a stationary environment and a stationary scenery. In order to
make sense of this, one needs to be able to compare the environment and
scenery at one state to those at another.
The simplest assumption is that there is a group $\Gamma$ of ``symmetries"
of the state space $\verts$ that acts transitively on $\verts$. Then
$\Gamma$ induces an action on functions on $\verts$, in particular, on
environments and sceneries. Restricting one's attention to the
$\sigma$-field $\inv$ of $\Gamma$-invariant events, one can ask whether
there is a stationary probability measure on $\inv$ that is equivalent to
the natural one.

In many cases of interest, there is such a stationary probability measure
that one can explicitly give. We present some general theorems of this
sort.  These are ``soft'' theorems, in contrast to most theorems in the
literature that describe more quantitative behavior of the processes.
There are some surprising phenomena even with such soft theorems. Compare
the following two examples:

\procl x.treebern
Consider a regular tree $T = (\verts, \edges)$ of degree 3 and fix $o \in
\verts$. Let $\Gamma$ be the group of automorphisms of $T$.
Declare each edge in $E$ ``open'' with probability 2/3 independently. Let
$\omega$ consist of the subgraph formed by the open edges.
Consider simple random walk starting at $o$ on the connected component
$C(o)$ of
$o$ in $\omega$. This has an equivalent stationary initial probability measure,
namely, the law of $\omega$ (product measure) biased by the degree of $o$
in $\omega$.
\endprocl

\procl x.treeend
With notation as above, let $\bdend$ be a fixed end of $T$.
Let $\Gamma_\bdend$ be the group of automorphisms on $T$ that fix $\bdend$.
This subgroup is also transitive on $\verts$.
However, in this case, simple random walk on $C(o)$ does not have any
stationary probability measure equivalent to the natural probability
measure: Let $Y(x)$ be the vertex in $C(x)$ that is closest to $\bdend$.
Let $w(n)$ denote the location of the walker at time $n$.
Let $\A_n$ be the event that $C(w(n))$ is infinite and $w(n) = Y(w(n))$; this
event is $\Gamma_\bdend$-invariant. Note that when the walker starts at
$o$, we have $C(w(n)) = C(o)$ and $Y(w(n)) = Y(o)$.
As time evolves, the probability of $\A_n$
tends to 0, yet the probability of $\A_0$ is positive.
\endprocl


It turns out that an important issue for finding a stationary measure is
whether $\Gamma$ is unimodular or not (see \ref s.def/ for the definition).
The group $\Gamma$ of \ref x.treebern/ is unimodular, but the group
$\Gamma_\bdend$ of \ref x.treeend/ is not.  In
many applications, $\verts$ is a countable group $\Gamma$ such as $\Z^d$,
in which case $\Gamma$ acts on itself by multiplication; since $\Gamma$ is
countable, it is unimodular.

In order to state one of our theorems, we need some notation.
The space of trajectories of the walk is $\verts^\N$. Let $(\Xi, \F)$ be a
measurable space which will be used to define the environment and the
scenery.

Define the shift $\SS:\verts^\N\to\verts^\N$ by
$$
(\SS w)(n):= w(n+1)
\,,
$$
and let 
$$
\SS(\xi, w):= (\xi, \SS w) \qquad\forall (\xi, w)\in\Xi\times \verts^\N
\,.
$$
For $\gamma\in\Gamma$, we set
$$
\gamma(\xi, w):= (\gamma \xi,\gamma w)
\,,
$$ 
where $(\gamma w)(n):=\gamma\bigl(w(n)\bigr)$.

A quadruple $(\Xi, \F, \P, \Gamma)$ is called a 
{\bf measure-preserving dynamical system} if
$\Gamma$ acts measurably on the measure space $(\Xi, \F, \P)$ preserving
the measure $\P$. We call a measurable function $p : \Xi \times \verts
\times \verts \to [0, 1]$, written $p : (\xi, x, y) \mapsto p_\xi(x, y)$, a
{\bf random environment} (from $\Xi$) if for all $\xi \in \Xi$ and all $x \in
\verts$, we have $\sum_{y \in \verts} p_\xi(x, y) = 1$.
The natural action of $\Gamma $ on $p$ is the one induced by the diagonal
one, $(\gamma p) (\xi, x, y) := p(\gamma^{-1} \xi, \gamma^{-1} x,
\gamma^{-1} y)$.
Unless otherwise stated, we shall use such actions implicitly.
Given $x \in \verts$ and a measurable map $\xi \mapsto \nu_\xi(x)$ from
$\Xi \to \CO{0, \infty}$,
let $\Ph_x$ denote the joint distribution on $\Xi \times \verts^\N$ of
$\xi$ biased by $\nu_\xi(x)$ and the trajectory of the Markov chain
determined by $p_\xi$ starting at $x$.
That is, if $\traj^x_\xi$ denotes the probability measure on $\verts^\N$
determined by $p_\xi$ with $w_0 = x$, then for all events $\ev A$, we have
$$
\Ph_x[\ev A]
:=
\int_{\Xi} d\P(\xi) \,\nu_\xi(x)\int_{(\xi, w) \in \ev A} d\traj^x_\xi(w) 
\,.
$$
Let $\inv$ be the $\sigma$-field of $\Gamma$-invariant events in $\Xi
\times \verts^\N$.
{\it We assume throughout this note that $\Gamma$ is a locally compact
group and that all stabilizers of elements of $\verts$ have finite Haar measure.}

The following theorem generalizes similar results in H\"aggstr\"om (1997),
H\"aggstr\"om and Peres (1999), Lyons and Peres (1998), and Lyons
and Schramm (1998).  

\procl t.main
Let $\verts$ be a countable set acted on by a transitive unimodular group
$\Gamma$.
\genassumps.
Suppose that $\nu : (\xi, x) \mapsto \nu_\xi(x)$ is a $\Gamma$-invariant
measurable mapping from $\Xi \times \verts \to \CO{0, \infty}$ such that for
each $\xi \in \Xi$, $\nu_\xi$ is a stationary distribution for the Markov
chain determined by $p_\xi$.
Then for any $o \in \verts$,
the restriction of $\Ph_o$ to the $\Gamma$-invariant $\sigma$-field is
an $\SS$-invariant measure; that is, $(\Xi \times \verts^\N, \inv, \Ph_o,
\SS)$ is a measure-preserving dynamical system.
If $\E[\nu_\buldot(o)] = 1$, then $\Ph_o$ is a probability measure.
\endprocl

As an example of an $\inv$-measurable function, we offer $p_\xi(w(0),
\cbuldot)$, the environment at the location of the walker. A function
$\scnry : \Xi \times \verts \to \R$ can be regarded as a random real-valued
scenery, where $\scnry(\xi, x)$ is the scenery at $x$ given by the outcome
$\xi$. If $\scnry$ is a $\Gamma$-invariant measurable function, then
$\scnry(\xi, w(0))$ is $\inv$-measurable. Thus, the theorem implies that
the walker will see a stationary scenery.


\procl x.alili
\procname{Alili (1994)} 
Let $\verts:= \Gamma := \Z$, $\Xi := (0, 1)^\Z$, 
$\P$ be any $\Z$-invariant measure on $\Xi$,
and for all $\xi\in\Xi$,
$$
p_\xi(x, y) := 
\cases{
\xi(x)   &if $y = x+1$,\cr
1 - \xi(x)  &if $y = x-1$,\cr
0   &otherwise.\cr
}
$$
Write $\rho(x) := \xi(x-1)/\xi(x)$.
Suppose that $A(x) := \sum_{n \ge x} \prod_{k=x+1}^n \rho(x) < \infty$ a.s.
Then $\nu_\xi(x) := (1+\rho(x)) A(x)$ is a stationary measure with $(\xi,
x) \mapsto \nu_\xi(x)$ being $\Z$-invariant.
\endprocl

\procl x.dsrw
Suppose that $G = (\verts, \edges)$ is a graph and $\Gamma$ is a closed
(vertex-)transitive group of automorphisms of $G$. Let $\P$ be a
$\Gamma$-invariant probability measure on $2^\edges$. That is, we choose a
random subgraph of $G$. 
The case that $\P$ is product measure, as in \ref x.treebern/,
is called {\bf Bernoulli
percolation}.
The random subgraph has connected components, often called ``percolation
clusters".
These clusters are of great interest. One method that has recently proven
quite powerful for studying the clusters is to use them for a random
environment (and/or scenery). Namely, let $D$ be the degree of vertices in
$G$. Denote the subgraph by $\omega$.
An RWRE called {\bf
delayed simple random walk} is defined via the transition probabilities
$p_\omega(x, y) := 1/D$ if $[x, y] \in \omega$ and $p_\omega(x, x) =
d_\omega(x)/D$, where $d_\omega(x)$ is the degree of $x$ in $\omega$.
This was introduced by H\"aggstr\"om (1997) and used also by H\"aggstr\"om
and Peres (1999), \BLSpert, and \LSindist.
If $\Gamma$ is unimodular, we take $\Xi := 2^\edges$ and $\nu \equiv 1$ in
\ref t.main/.
\endprocl

\procl x.srw
In the same setting as \ref x.dsrw/, consider the transition probabilities
$p_\omega(x, y) := 1/d_\omega(x)$ if $[x, y] \in \omega$ and $d_\omega(x)
\ne 0$, with $p_\omega(x, x) = 1$ if $d_\omega(x) = 0$. This is called {\bf
simple random walk on percolation clusters}.
In this case, we take $\nu_\omega(x) := d_\omega(x)$ if $d_\omega(x)
\ne 0$ and $\nu_\omega(x) := 1$ if $d_\omega(x) = 0$.
The paper by
\BLSpert\ gives a number of potential-theoretic properties of simple random
walk on percolation clusters.
\endprocl

\bsection{Definitions}{s.def}

Let $\verts$ be a countable set.
If $\Gamma$ acts on $\verts$ (on the left),
we say that $\Gamma$ is {\bf
transitive} if for every $x, y\in\verts$, there is a $\gamma\in\Gamma$
with $\gamma x=y$.  
If the orbit space $\Gamma\backslash\verts$ is finite,
then $\Gamma$ is {\bf quasi-transitive}.

Recall that on every locally compact group $\Gamma$, there is
a unique (up to a constant scaling factor) Borel measure $|\,\cbuldot\,|$
that, for every $\gamma\in\Gamma$, is invariant under left multiplication
by $\gamma$; this measure is called (left) {\bf Haar measure}.
The group is {\bf unimodular} if Haar measure is also invariant
under right multiplication.
For example, when $\Gamma$ is countable, the Haar measure is (a constant
times) counting measure, so $\Gamma$ is unimodular.
Let 
$$
S(x) := \{\gamma \in \Gamma \st \gamma x = x\}
$$
denote the stabilizer of $x$.
We shall write 
$$
m(x) := |S(x)|
\,.
$$
It is not hard to show that a group
$\Gamma$ with stabilizers of finite Haar measure
is unimodular iff $m(\cbuldot)$ is $\Gamma$-invariant iff
for all $x$ and $y$ in the same orbit,
$
|S(x)y| = |S(y)x| 
$
(see Trofimov (1985)).

H\"aggstr\"om (1997) introduced the Mass-Transport Principle
in studying percolation on regular trees.
Following is a generalization.

\procl l.qtrNonuni
Let $\Gamma$ act quasi-transitively on $\verts$ and $f : \verts \times
\verts \to [0, \infty]$ be invariant under the diagonal action of $\Gamma$.
Choose a complete set $\{o_1, \ldots, o_L\}$ of representatives in
$\verts$ of the orbits of $\Gamma$ and write $m_i := m(o_i)$. Then
$$
\sum_{i=1}^L \sum_{z \in \verts} f(o_i, z)
= \sum_{j=1}^L 1/m_j \sum_{y \in \verts} f(y, o_j) m(y)\,.
$$
\endprocl

\noindent See Cor.~3.7 of \BLPSgip.

\bsection{Proofs}{s.proof}

\ref t.main/ generalizes as follows to quasi-transitive actions:

\procl t.qtr
Let $\verts$ be a countable set acted on by a quasi-transitive unimodular
group $\Gamma$.
Let $\{o_1, \dots, o_L\}$ be a complete set of representatives of $\Gamma
\backslash \verts$ and write $m_i := m(o_i)$.
\genassumps.
Suppose that $\nu : (\xi, x) \mapsto \nu_\xi(x)$ is a $\Gamma$-invariant
measurable mapping from $\Xi \times \verts \to \CO{0, \infty}$ such that for
each $\xi \in \Xi$, $\nu_\xi$ is a stationary distribution for the Markov
chain determined by $p_\xi$.
Write
$$
\Ph := \sum_{i=1}^L m_i^{-1} \Ph_{o_i}
\,.
$$
Then the restriction of $\Ph$ to the $\Gamma$-invariant $\sigma$-field is an
$\SS$-invariant measure.
If 
$$
\sum_i m_i^{-1} \E[\nu_\buldot(o_i)] = 1
\,,
$$
then $\Ph$ is a probability measure.
\endprocl

Still more generally, we may remove the hypothesis that $\Gamma$ be
unimodular by means of the following modification: 

\procl t.gen
Let $\verts$ be a countable set acted on by a quasi-transitive 
group $\Gamma$.
Let $\{o_1, \dots, o_L\}$ be a complete set of representatives of $\Gamma
\backslash \verts$ and write $m_i := m(o_i)$.
\genassumps.
Suppose that $\nu : (\xi, x) \mapsto \nu_\xi(x)$ is a $\Gamma$-invariant
measurable mapping from $\Xi \times \verts \to \CO{0, \infty}$ such that 
for each $\xi \in \Xi$, $x \mapsto m(x) \nu_\xi(x)$ is a stationary
distribution for the Markov chain determined by $p_\xi$.
Write
$$
\Ph := \sum_{i=1}^L \Ph_{o_i}
\,.
$$
Then the restriction of $\Ph$ to the $\Gamma$-invariant $\sigma$-field is an
$\SS$-invariant measure.
If 
$$
\sum_i \E[\nu_\buldot(o_i)] = 1
\,,
$$
then $\Ph$ is a probability measure.
\endprocl

Note that this incorporates \ref t.qtr/ because when $\Gamma$ is
unimodular, the function $(\xi, x) \mapsto m(x) \nu_\xi(x)$ is
$\Gamma$-invariant.

\proof
Let $F$ be a $\Gamma$-invariant function on $\Xi \times \verts^\N$. We must
show that
$\int d\Ph \,F \circ \SS = \int d\Ph \,F$.

Set
$$
f(x, y; \xi) := \nu_\xi(x) p_\xi(x, y) \int d\traj_\xi^y(w) \,F(\xi, w)
\,.
$$
Thus, we have
$$
\int d\Ph \,F \circ \SS 
=
\sum_{i=1}^L \sum_{y \in \verts} \int d\P(\xi) \, f(o_i, y; \xi)
\,.
$$
Our assumptions imply that $f$, and hence $\E[f(x, y; \cbuldot)]$,
is $\Gamma$-invariant. 
Consequently, \ref l.qtrNonuni/ gives
\begineqalno
\int d\Ph \,F \circ \SS 
&=
\sum_{j=1}^L \sum_{y \in \verts} \int
d\P(\xi) \, m(y) f(y, o_j; \xi) /m_j
\cr&=
\sum_{j=1}^L \int d\P(\xi) \sum_{y \in \verts} 
\nu_\xi(y) m(y) p_\xi(y, o_j) /m_j \int d\traj_\xi^{o_j}(w) \,F(\xi, w)
\cr&=
\sum_{j=1}^L \int d\P(\xi) \, \nu_{\xi}(o_j)
\int d\traj_\xi^{o_j}(w) \,F(\xi, w)
=
\int d\Ph \,F
\,.
\Qed
\endeqalno

\procl x.brw
Suppose that $G = (\verts, \edges)$ is a graph and $\Gamma$ is a closed
quasi-transitive group of automorphisms of $G$. Let $\P$ be a
$\Gamma$-invariant probability measure on $2^\edges$. 
Write
$$
\alpha(x) := \sum_{[x, y] \in \edges} \sqrt{m(y)/m(x)}
\,.
$$
Given $\omega \in 2^\edges$,
consider the transition probabilities $p_\omega(x, y) := \alpha(x)^{-1}
\sqrt{m(y)/m(x)}$ for $[x, y] \in \omega$ and $p_\omega(x, x) := 1 -
\sum_{[x, y] \in \omega} p_\omega(x, y)$. 
The resulting Markov chain on $\omega$ is reversible with stationary
measure $x \mapsto m(x)\alpha(x)$.
In the unimodular transitive case, this Markov chain is delayed simple
random walk.
Whether $\Gamma$ is unimodular or not, we may take $\Xi := 2^\edges$ and
$\nu_\omega(x) := \alpha(x)$ in \ref t.gen/.
\endprocl

\noindent{\bf Acknowledgement.} We are grateful to Yuval Peres
for fruitful conversations.

\beginreferences

Alili, S. (1994) Comportement asymptotique d'une marche al\'eatoire en
environnement al\'eatoire, {\it C. R. Acad. Sci. Paris S\'er. I Math.} {\bf
319}, 1207--1212.


Benjamini, I., Lyons, R., Peres, Y., \and Schramm, O. (1999) 
 Group-invariant percolation on graphs,
{\it Geom. Funct. Anal.} {\bf 9}, 29--66.

Benjamini, I., Lyons, R. \and Schramm, O. (1998) Percolation perturbations
in potential theory and random walks, in {\it Random Walks and Discrete
Potential Theory (Cortona, 1997), Sympos. Math.}, Cambridge
Univ. Press, Cambridge, 1998, to appear.

H\"aggstr\"om, O. (1997) Infinite clusters in dependent automorphism
invariant percolation on trees, {\it Ann. Probab.} {\bf 25}, 1423--1436.

H\"aggstr\"om, O. \and Peres, Y. (1999) Monotonicity of uniqueness for
percolation on Cayley graphs: all infinite clusters are born
simultaneously, 
{\it Probab. Theory Rel. Fields} {\bf 113}, 273--285.

Lang, R. \and Nguyen, X.-X. (1983) Strongly correlated random fields as
observed by a random walker, {\it Z. Wahrsch. Verw. Gebiete} {\bf 64},
327--340.

Lyons, R. \and Peres, Y. (1998) {\it Probability on Trees and Networks}.
Cambridge University Press, in preparation. Current version available 
at \hfill\break {\tt http://php.indiana.edu/\~{}rdlyons/}.

Lyons, R. \and Schramm, O. (1998) Indistinguishability of percolation
clusters, {\it Ann. Probab.}, to appear.

Solomon, F. (1975) Random walks in a random environment, {\it Ann.
Probab.} {\bf 3}, 1--31.

\endreferences

\filbreak
\begingroup
\eightpoint\sc
\parindent=0pt\baselineskip=10pt

\def\emailwww#1#2{\par\qquad {\tt #1}\par\qquad {\tt #2}\smallskip}

Department of Mathematics,
Indiana University,
Bloomington, IN 47405-5701, USA
\emailwww{rdlyons@indiana.edu}
{http://php.indiana.edu/\~{}rdlyons/}

Mathematics Department,
The Weizmann Institute of Science,
Rehovot 76100, Israel
\emailwww{schramm@wisdom.weizmann.ac.il}
{http://www.wisdom.weizmann.ac.il/\~{}schramm/}

\endgroup

\bye